\nonstopmode
\nonstopmode
\documentclass[]{amsart}
\usepackage{verbatim}
\usepackage{amssymb,bbm}
\usepackage{graphicx}
\usepackage{psfrag}
\usepackage{ifpdf}

\ifpdf
\usepackage[notightpage]{pst-pdf}
\fi

\allowdisplaybreaks
\usepackage[all]{xy}
\newdir{ >}{{}*!/-10pt/@{>}}
\newdir^{ (}{{}*!/-5pt/@^{(}}
\newdir_{ (}{{}*!/-5pt/@_{(}}
\def\cgaps#1{}
\def\Cgaps#1{}
\def\AMStag#1{}
\def\AMSunderset#1\to#2{\underset{#1}{#2}}
\def\AMSoverset#1\to#2{\overset{#1}{#2}}
\def\undersetbrace#1\to#2{\underbrace{#2}_{#1}}
\def\oversetbrace#1\to#2{\overbrace{#2}^{#1}}
\def\3{\ss}

\newcommand{\nmb}[2]{\ifx!#1{\ref{nmb:#2}}%
\else\if.#1{\label{nmb:#2}}%
\else\if0#1{\label{nmb:#2}}%
\else{{#2}}%
\fi\fi\fi}
\swapnumbers
\newtheorem{proposition}[subsection]{Proposition}
\newtheorem*{proposition*}{Proposition}
\newtheorem{theorem}[subsection]{Theorem}
\newtheorem*{theorem*}{Theorem}
\newtheorem{lemma}[subsection]{Lemma}
\newtheorem*{lemma*}{Lemma}

\newtheorem*{corollary*}{Corollary}

\newtheorem*{result*}{Result}
\newenvironment{demo}[1]{\par\smallskip\noindent{\bf #1.}}{\par\smallskip}

\def\East#1#2{-\raisebox{0.1pt}{$\mkern-16mu\frac{\;\;#1\;}{\;\;#2\;}\mkern-16mu$}\to}

\def\cit#1#2{\ifx#1!\cite{2}\else#2\fi} 
\def\idx{}               
\def\ign#1{}             
\def\o{\circ} 
\def\X{\mathfrak X} 
\def\al{\alpha} 
\def\be{\beta} 
 
\def\de{\delta} 
\def\ep{\varepsilon}

\def\la{\lambda}

\def\ph{\varphi} 
 
\def\ps{\psi} 
\def\om{\omega}

\def\i{^{-1}} 
\def\x{\times} 
\def\g{{\mathfrak g}} 
\def\p{\partial}
\def\ad{\operatorname{ad}}

\let\on=\operatorname
\def\supp{{\on{supp}}}

\DeclareMathOperator{\Id}{Id}
\DeclareMathOperator{\Diff}{Diff}

\newcommand{\R}{\mathbb{R}}

\newcommand{\wt}[1]{\widetilde{#1}}

\newcommand{\ol}[1]{\overline{#1}}

\newcommand{\ud}{\,\mathrm{d}}

\newcommand{\one}{\mathbbm{1}}

\title[Vanishing geodesic distance for KdV]
{Vanishing geodesic distance for the Riemannian metric with geodesic equation the KdV-equation}
\author{Martin Bauer, Martins Bruveris, Philipp Harms, Peter W. Michor}
 
\address{
Martin Bauer, Philipp Harms, Peter W.\ Michor: Fakult\"at f\"ur Mathematik, Universit\"at Wien, 
Nordbergstrasse 15, A-1090 Wien, Austria.\newline\indent
Martins Bruveris: Dep. of Mathematics, Imperial College, London SW7~2AZ, UK.
}
\email{bauer.martin@univie.ac.at}
\email{m.bruveris08@imperial.ac.uk}
\email{philipp.harms@univie.ac.at}
\email{peter.michor@esi.ac.at }
\date{{\today} } 
\thanks{All authors were supported by `Fonds zur
F\"orderung der wissenschaftlichen                    
Forschung, Projekt P~21030'}
\begin{document}
\begin{abstract}  
The Virasoro-Bott group endowed with the right-invariant $L^2$-metric (which is a weak Riemannian metric) 
has the KdV-equation as geodesic equation. We prove that this metric space has vanishing geodesic distance.
\end{abstract} 

\maketitle

\section{\nmb0{1} Introduction}

It was found in \cite{OvsienkoKhesin87} that  a curve in the Virasoro-Bott group is a geodesic for 
the right invariant $L^2$-metric if and only if its right logarithmic derivative is a 
solution of the Korteweg-de Vries equation, see \nmb!{2.3}. 
Vanishing geodesic distance for weak Riemannian metrics on infinite dimensional manifolds was first 
noticed on shape space $\on{Imm}(S^1,\mathbb R^2)/\on{Diff}(S^1)$ for the $L^2$-metric in 
\cite[3.10]{MM06}.  
In \cite{Michor2005} this result was shown to hold for the general shape space $\on{Imm}(M,N)/\on{Diff}(M)$ for any compact 
manifold $M$ and Riemannian manifold $N$, and also for the right invariant $L^2$-metric on 
each full diffeomorphism group with compact 
support $\on{Diff}_c(N)$. In particular, Burgers' equation is related to the geodesic equation of the 
right invariant $L^2$-metric on $\on{Diff}(S^1)$ or $\on{Diff}_c(\mathbb R)$ and it thus also has 
vanishing geodesic distance. We even have 

\begin{result*} \cite{Michor2005} The weak Riemannian $L^2$-metric on each connected component 
of the total space $\on{Imm}(M,N)$ for 
a compact manifold $M$ and a Riemannian manifold $(N,g)$ has vanishing geodesic distance.
\end{result*}

This result is not spelled out in \cite{Michor2005} but it follows from there: Given two immersions 
$f_0,f_1$ in the same connected component, 
we first connect their shapes $f_0(M)$ and $f_1(M)$ by a curve of length $<\ep$ in the shape 
space $\on{Imm}(M,N)/\on{Diff}(M)$ and take the horizontal lift to get a curve of length $<\ep$ 
from $f_0$ to an immersion $f_1\o \ph$ in the connected component of the orbit through $f_1$. 
Now we use the induced metric 
$f_1^*g$ on $M$ and the right invariant $L^2$-metric induced on $\on{Diff}(M)_0$ to get a curve in 
$\on{Diff}(M)$ of length $<\ep$ connecting $\ph$ with $\on{Id}_M$. Evaluating at $f_1$ we get curve 
in $\on{Imm}(M,N)$ of length $<\ep$ connecting $f_1\o\ph$ with $f_1$.

In this article we show that the right invariant $L^2$-metric on the Virasoro-Bott groups (see 
\nmb!{2.1}) has vanishing geodesic distance. This might be related to the fact that the Riemannian 
exponential mapping is not a diffeomorphism near $0$, see \cite{Constantin2003} for $\on{Diff}(S^1)$ 
and \cite{CKKT07} for the Virasoro group over $S^1$. 
See \cite{Misiolek97} for information on conjugate points along geodesics.

\section{\nmb0{2} The Virasoro-Bott groups} 
 
\subsection{\nmb.{2.1} The Virasoro-Bott groups  } 
Let $\on{Diff}_{\mathcal{S}}(\mathbb R)$ be the group of diffeomorphisms 
of $\mathbb R$ which rapidly fall to the identity.
This is a regular Lie group, see \cite[6.4]{MGeomEvol06}.
The mapping 
\begin{gather*} 
c:\on{Diff}_{\mathcal{S}}(\mathbb R)\x \on{Diff}_{\mathcal{S}}(\mathbb R)\to \mathbb R
\\
c(\ph,\ps):=\frac12\int\log(\ph\o\ps)'d\log\ps'  
     = \frac12\int\log(\ph'\o\ps)d\log\ps' 
\end{gather*}
satisfies $c(\ph,\ph\i)=0$, $c(\on{Id},\ps)=0$, $c(\ph,\on{Id})=0$ and is a 
smooth group cocycle, called the Bott cocycle: 
\begin{displaymath}
c(\ph_2,\ph_3) - c(\ph_1\o\ph_2,\ph_3)+c(\ph_1,\ph_2\o\ph_3) -
c(\ph_1,\ph_2) =0.
\end{displaymath}

The corresponding central extension group
$\on{Vir}:=\mathbb R\x_c\on{Diff}_{\mathcal S}(\mathbb R)$, called the Virasoro-Bott group, is
a trivial $\mathbb R$-bundle  
$\mathbb R\x\on{Diff}_{\mathcal S}(\mathbb R)$ that becomes a regular Lie  
group relative to the operations  
\begin{displaymath}
\binom{\ph}{\al}\binom{\ps}{\be}  
     =\binom{\ph\o\ps}{\al+\be+c(\ph,\ps)},\quad 
\binom{\ph}{\al}\i=\binom{\ph\i}{-\al}\quad 
\ph, \ps \in \on{Diff}_{\mathcal S}(\mathbb R),\; \al,\be \in \mathbb R . 
\end{displaymath}
Other versions of the Virasoro-Bott group are the following:
$\mathbb R\x_c \on{Diff}_c(\mathbb R)$ where $\on{Diff}_c(\mathbb R)$ is the group of all 
diffeomorphisms with compact support, or the periodic case $\mathbb R\x_c\on{Diff}^+(S^1)$.
One can also apply  the homomorphism $\exp(i\al)$ to the center and replace it by $S^1$. 
To be specific, we shall treat the most difficult case $\on{Diff}_{\mathcal{S}}(\mathbb R)$ 
in this article. All other 
cases require only obvious minor changes in the proofs.

\subsection{\nmb.{2.2} The Virasoro Lie algebra }
The Lie algebra of the Virasoro-Bott
group $\mathbb R\x_c\on{Diff}_{\mathcal S}(\mathbb R)$
is $\mathbb R\x \X_{\mathcal S}(\mathbb R)$ (where 
$\X_{\mathcal S}(\mathbb R)=\mathcal S(\mathbb R)\p_x$) 
with the Lie bracket 
\begin{equation*}
\left[\binom{X}{a},\binom{Y}{b} \right]
=\binom{-[X,Y]}{\om(X,Y)} = \binom{X'Y-XY'}{\om(X,Y)}
\end{equation*}
where
\begin{displaymath}
\om(X,Y)=\om(X)Y=\int X'dY'=\int X'Y''dx =  
\tfrac12\int \det\begin{pmatrix} X'& Y'\\ X''&Y''\end{pmatrix}\,dx, 
\end{displaymath}
is the \idx{\it Gelfand-Fuks Lie algebra cocycle}  
$\om:\g\x \g\to \mathbb R$, which is a bounded skew-symmetric bilinear mapping
satisfying the cocycle condition
\begin{displaymath}
\om([X,Y],Z)+\om([Y,Z],X)+\om([Z,X],Y)=0.
\end{displaymath}
It is a generator of the 1-dimensional bounded Chevalley cohomology  
$H^2(\g,\mathbb R)$ for any of the Lie algebras 
$\g=\X(\mathbb R)$, $\X_c(\mathbb R)$, or $\X_{\mathcal S}(\mathbb R)=\mathcal{S}(\mathbb R)\p_x$.
The Lie algebra of the Virasoro-Bott Lie group is thus the central extension 
$\mathbb R\x_\om \X_{\mathcal S}(\mathbb R)$ induced by this cocycle.
We have $H^2(\X_c(M),\mathbb R)=0$ for each
finite dimensional manifold of dimension $\ge 2$ (see \cite{Fuks86}),  
which blocks the way to  
find a higher dimensional analog of the Korteweg-de Vries  
equation in a way similar to that sketched below. 

To complete the description, we add the adjoint action: 
\begin{equation*}
\on{Ad}\binom{\ph}{\al}\,\binom{Y}{b} 
=  \binom{\on{Ad}(\ph)Y=\ph_*Y=T\ph\o Y\o \ph\i}
{b+\int S(\ph)Y\,dx}
\end{equation*}
where the \idx{\it Schwartzian
derivative} $S$ is given by
\begin{align*}
S(\ph) &
= \Bigl(\frac{\ph''}{\ph'}\Bigr)' -\frac{1}2\Bigl(\frac{\ph''}{\ph'}\Bigr)^2
= \frac{\ph'''}{\ph'} -\frac{3}2\Bigl(\frac{\ph''}{\ph'}\Bigr)^2
= \log(\ph')'' -\frac12(\log(\ph')')^2
\end{align*}
which measures the deviation of $\ph$ from being a M\"obius transformation:
\begin{displaymath}
S(\ph)=0 \iff \ph(x)=\frac{ax+b}{cx+d}\text{  for }
\begin{pmatrix} a & b \\ c & d\end{pmatrix} \in SL(2,\mathbb R).
\end{displaymath}
The Schwartzian derivative of a
composition and an inverse follow from the action property:
\begin{displaymath}
S(\ph\o \ps) = (S(\ph)\o\ps)(\ps')^2 + S(\ps),\quad
S(\ph\i)=-\frac{S(\ph)}{(\ph')^2}\o\ph\i
\end{displaymath}

\subsection{\nmb.{2.3} The right invariant $L^2$-metric and the KdV-equation}
We shall use the $L^2$-inner product on $\mathbb R\x_\om \X_{\mathcal{S}}(\mathbb R)$: 
\begin{equation*}
\left\langle \binom{X}{a},\binom{Y}{b}\right\rangle :=  
\int XY\,dx + ab. 
\end{equation*}
We use the induced right invariant weak Riemannian metric on the Virasoro group. 

According to \cite{Arnold66}, see \cite{MR98} for a proof in the notation and setup used here,
a curve 
$t \mapsto \binom{\ph(t,\quad)}{\al(t)}$
in the Virasoro-Bott group is a geodesic if and only if
\begin{align*}
\binom{u_t}{a_t}&=-\ad\binom{u}{a}^\top\binom{u}{a} 
     =\binom{-3u_xu-au_{xxx}}{0}\quad\text{  where } 
\\
\binom{u(t)}{a(t)}&=\p_s\binom{\ph(s)}{\al(s)}.\binom{\ph(t)\i}{-\al(t)}\Bigr|_{s=t}
=\p_s\binom{\ph(s)\o\ph(t)\i}{\al(s)-\al(t)+c(\ph(s),\ph(t)\i)}\Bigr|_{s=t},
\\
\binom{u}{a}
&=\binom{\ph_t\o\ph\i}{\al_t-\int \frac{\ph_{tx}\ph_{xx}}{2\ph_x^2}dx},
\end{align*}
since we have
\begin{align*}
2\p_s &c(\ph(s),\ph(t)\i)|_{s=t} = \p_s\int\log(\ph(s)'\o
  \ph(t)\i)\,d\log((\ph(t)\i)')|_{s=t}
\\&
= \int \frac{\ph_{t}(t)'\o \ph(t)\i}{\ph(t)'\o
\ph(t)\i}\left(-\frac{\ph(t)''\o\ph(t)\i}{(\ph(t)'\o\ph(t)\i)^2} \right)\,dx
\\&
= -\int \Big(\frac{\ph_{t}'\ph''}{(\ph')^2}\Big)(t)dy
= -\int \Big(\frac{\ph_{tx}\ph_{xx}}{\ph_x^2}\Big)(t)dx.
\end{align*}

Thus $a$ is a constant in time and the geodesic equation  
is hence the {\it Korteweg-\-de~Vries equation}
\begin{equation*}
u_t+3u_xu+au_{xxx}=0
\end{equation*}
with its natural companions
\begin{displaymath}
\ph_t=u\o\ph,\qquad \al_t = a+\int \frac{\ph_{tx}\ph_{xx}}{2\ph_x^2}dx. 
\end{displaymath}

To be complete, we add the invariant momentum mapping $J$ 
with values in the Virasoro algebra (via the weak Riemannian 
metric). We
need the transpose of the adjoint action:
\begin{align*}
&\left\langle \on{Ad}\binom{\ph}{\al}^\top\,\binom{Y}{b},\binom{Z}{c} \right\rangle
=\left\langle \binom{Y}{b},\on{Ad}\binom{\ph}{\al}\,\binom{Z}{c} \right\rangle
\\&\qquad
=\left\langle \binom{Y}{b},\binom{\ph_*Z}{c+\int S(\ph)Z\,dx}\right\rangle
\\&\qquad
=\int Y((\ph'\o\ph\i)(Z\o\ph\i)\,dx +bc+\int bS(\ph)Z\,dx
\\&\qquad
=\int ((Y\o\ph)(\ph')^2 +bS(\ph))Z\,dx + bc
\end{align*}
Thus, the invariant momentum mapping
is given by  
\begin{equation*}
J\left(\binom{\ph}{\al},\binom{Y}{b} \right) =
\on{Ad}\binom{\ph}{\al}^\top\binom{Y}{b} 
=\binom{(Y\o\ph)(\ph')^2 +bS(\ph)}{b}.
\end{equation*}
Along a geodesic $t\mapsto g(t,\quad)=\binom{\ph(t,\quad)}{\al(t)}$, 
the momentum 
\begin{equation*}
J\left(\binom{\ph}{\al},\binom{u=\ph_t\o \ph\i}{a}\right) 
= \binom{(u\o \ph)\ph_x^2 +aS(\ph)}{a} = \binom{\ph_t\ph_x^2+aS(\ph)}{a}
\end{equation*}
is constant in $t$.

\subsection{\nmb.{2.4} Lifting curves to the Virasoro-Bott group}
We consider the extension 
$$
\mathbb R \East{i}{} \mathbb R\x_c \on{Diff}_{\mathcal S}(\mathbb R) \East{p}{} 
\on{Diff}_{\mathcal S}(\mathbb R).
$$
Then $p$ is a Riemannian submersion for the right invariant $L^2$-metric on 
$\on{Diff}_{\mathcal S}(\mathbb R)$, i.e., $Tp$ is an isometry on the orthogonal complements 
of the fibers. These complements are not integrable; in fact, the curvature of the corresponding 
principal connection is given by the Gelfand-Fuks cocycle. 
For any curve $\ph(t)$ in $\on{Diff}_{\mathcal S}(\mathbb R)$ its horizontal lift is given by 
$$ \binom{\ph(t)}{a(t)           =           a(0)          -\int_0^t
\int\frac{\ph_{tx}\ph_{xx}}{\ph_x^2} \,d x \,d t}
$$
since the right translation to $(\on{Id},0)$ of its velocity should have zero vertical component, see
\nmb!{2.3}.
The horizontal lift has the same length and energy as $\ph$.

\section{\nmb0{3} Vanishing of the geodesic distance}

\begin{theorem}\nmb.{3.1}
On all Virasoro-Bott groups mentioned in \nmb!{2.1} geodesic distance for the right invariant $L^2$-metric vanishes.
\end{theorem}

The rest of this section is devoted to the proof of theorem \nmb!{3.1} for the most difficult case 
$\mathbb R\x_c\on{Diff}_{\mathcal S}(\mathbb R)$.

\begin{proposition}\nmb.{3.2}
Any two diffeomorphisms in $\Diff_{\mathcal S}(\R)$ can be connected by 
a path with arbitrarily short length for the right invariant $L^2$-metric.
\end{proposition}

In \cite{Michor2005} for $\on{Diff}_c(\mathbb R)$ it was first shown 
that there exists one non-trivial diffeomorphism which can 
be connected to $\Id$ with arbitrarily small length. 
Then, it was shown that the diffeomorphisms with this 
property form a normal subgroup. Since  $\on{Diff}_c(\mathbb R)$ is a simple group this concluded 
the proof.
But $\on{Diff}_{\mathcal S}(\mathbb R)$ is not a simple group since $\on{Diff}_c(\mathbb R)$ 
is a normal subgroup. So, we have to elaborate on the proof of \cite{Michor2005} as follows. 

\begin{demo}{Proof}
We show that any rapidly decreasing diffeomorphism can be connected
to the identity by an arbitrarily short path.
We will write this diffeomorphism as $\Id+g$, where $g \in \mathcal S(\R)$ is a rapidly 
decreasing function with $g'>-1$. 
For $\lambda=1-\ep <1$ we define
\[ \ph(t,x) = x + \max(0, \min(t - \lambda x, g(x))) 
    - \max(0,\min(t+\la x,-g(x))). \]
This is a (non-smooth) path defined for $t \in (-\infty,\infty)$ 
connecting the identity in $\on{Diff}_{\mathcal S}(\mathbb R)$ with the diffeomorphism $(\Id+g)$.
We define $\ps(t,x)=\ph(\tan(t),x) \star G_{\ep}(t,x)$, 
where $G_\ep(t,x) = \frac{1}{\ep^2} G_1(\frac{t}{\ep},\frac{x}{\ep})$
is a smoothing kernel with $\operatorname{supp}(G_\ep) \subseteq B_{\ep}(0)$ 
and $\iint G_\ep \ud x \ud t = 1$. 
Thus $\ps$ is a smooth path defined on the finite interval $-\tfrac{\pi}2<t<\tfrac{\pi}2$
connecting the identity in $\Diff_S(\mathbb R)$ with a diffeomorphism 
arbitrarily close to $(\Id+g)$ for $\ep$ small.
(Compare figure~\ref{fig} for an illustration.) 

The $L^2$-energy of $\ps$ is
\begin{align*}
E(\ps)&=\int_{-\tfrac{\pi}2}^{\tfrac{\pi}2} \int_{\R} (\ps_t \o \ps\i)^2 \ud x \ud t = 
\int_{-\tfrac{\pi}2}^{\tfrac{\pi}2} \int_{\R} \ps_t^2 \ps_x \ud x \ud t 
\end{align*}
where $\ps\i(t,x)$ stands for $\ps(t,\quad)\i(x)$.
We have
\begin{align*}
\p_a \max(0, \min(a, b)) &= \one_{0 \leq a \leq b},&
\p_b \max(0, \min(a, b)) &= \one_{0 \leq b \leq a},
\end{align*}
and therefore
\begin{align*}
\ps_x(t,x) &= \ph_x(\tan(t),x) \star G_{\ep} \\&=
	\big(1-\la \one_{0\leq \tan(t)-\la x \leq g(x)}\;\,\,
	+g'(x)\one_{0\leq g(x) \leq \tan(t)-\la x} \\&\qquad\;
	-\la \one_{0\leq \tan(t)+\la x \leq -g(x)}
	+g'(x)\one_{0\leq -g(x) \leq \tan(t)+\la x}\big)\star G_{\ep},\\
\ps_t(t,x) &= \big((1+\tan(t)^2) \ph_t(\tan(t),x)\big) \star G_{\ep} \\&=
	\big((1+\tan(t)^2) (\one_{0\leq \tan(t)-\la x \leq g(x)}
	-\one_{0\leq \tan(t)+\la x \leq -g(x)})\big) \star G_{\ep}.
\end{align*}
Note that these functions have disjoint support when $\ep=0, \la=1-\ep=1$. 

\noindent {\bf Claim.} The mappings $\ep \mapsto \ps_t$ and $\ep \mapsto (\ps_x-1)$ are 
continuous into each $L^p$ with $p$ even. 
(The proofs are simpler when $p$ is even because there are no absolute values to be taken care of.)
To prove the claim, we calculate
\begin{align*}&
\int_{-\tfrac{\pi}2}^{\tfrac{\pi}2}\int_{\R} \big((1+\tan(t)^2) \ph_t(\tan(t),x)\big)^p \ud x \ud t 
=
\iint_{\R^2} \ph_t(t,x)^p (1+t^2)^{p-1} \ud x \ud t 
\\&\qquad=
\iint_{\R^2} (\one_{0\leq t-\la x \leq g(x)}+\one_{0\leq t+\la x \leq -g(x)}) (1+t^2)^{p-1} \ud x \ud t 
\\&\qquad=
\int_{g(x) \geq 0}\int_{\la x}^{\la x + g(x)} (1+t^2)^{p-1} \ud t \ud x
+ \int_{g(x) < 0} \int_{\la x + g(x)}^{\la x} (1+t^2)^{p-1} \ud t \ud x
\\&\qquad=
\int_{\R} \left| F(t)\big|_{t=\la x}^{t=\la x+g(x)} \right| \ud x
= \int_{\mathbb R}|F(\la x + g(x)) - F(\la x)| \ud x,
\end{align*}
where $F(\la x + g(x))- F(\la x)$ is a polynomial without constant term in $g(x)$ with coefficients 
also powers of $\la x$. 
Integrals of the form
$\int_{\mathbb R}|(\la x)^{k_1} g(x)^{k_2}| \ud x$
with $k_1 \geq 0, k_2 >0$ are finite and continuous in $\la=1-\ep$ since $g$ is rapidly decreasing.
This shows that $\| (1+\tan(t)^2) \ph_t(\tan(t),x) \|_p$ 
depends continuously on $\ep$ .
Furthermore the sequence $(1+\tan(t)^2) \ph_t(\tan(t),x)$ converges almost everywhere for $\ep \to 0$, thus
it also converges in measure. By the theorem of Vitali, 
this implies convergence in $L^p$, see for example \cite[theorem~16.6]{Schilling2006}. 
Convolution with $G_{\ep}$ acts as approximate unit in each $L^p$, which proves the claim for $\ps_t$. 
For $\ps_x-1$ it follows similarly. 

The above claim implies that
\begin{align*}
E(\ps)&=\int_{-\tfrac{\pi}2}^{\tfrac{\pi}2} \int_{\R} \ps_t^2 \ps_x \ud x \ud t=
\int_{-\tfrac{\pi}2}^{\tfrac{\pi}2} \int_{\R} \ps_t^2 (\ps_x-1) \ud x \ud t
+\int_{-\tfrac{\pi}2}^{\tfrac{\pi}2} \int_{\R} \ps_t^2 \ud x \ud t
\end{align*}
viewed as a mapping on $L^4\x L^4 \x L^2$ (first summand) and 
on $L^2 \x L^2$ (second summand) is continuous in $\ep$. It also vanishes at $\ep=0$ since then $\ps_x$ 
and $\ps_t$ have disjoint support. 
The Cauchy-Schwarz inequality $L(\ps)^2 < \pi E(\ps)$ implies that $L(\ps)$ goes 
to zero as well.
Ultimately,
$\ps(\tfrac{\pi}2)=(\Id+g)\star G_{\ep}$ is arbitrarily
close to $\Id+g$.
\qed\end{demo}

\begin{lemma}\nmb.{3.3}
For any $a\in \R$ there exists an arbitrarily short path connecting $\binom{\Id}{0}$ and 
$\binom{\Id}{a}$, i.e., 
$\on{dist}_{\on{Vir}}^{L^2}\big(\binom{\Id}{0},\binom{\Id}{a}\big)=0.$ 
\end{lemma}
\begin{demo}{Proof}
The aim of the following argument is to construct a family of paths 
in the diffeomorphism group, parametrized by $\ep$, with the following
properties: 
all paths in the family start and end at the identity and their length in the
diffeomorphism group with respect to the $L^2$ metric tends to $0$ as $\ep
\to 0$. By letting $\ep$ be time-dependent, we are able to control the endpoint $a(T)$ of
the horizontal lift for certain diffeomorphisms.

\begin{figure}[h]
\begin{psfrags}
\providecommand{\PFGstyle}{}%
\psfrag{Idg}[cc][cc]{\PFGstyle $\operatorname{Id}+g$}%
\psfrag{Id}[tl][tl]{\PFGstyle $\operatorname{Id}$}%
\psfrag{Phit}[cr][cr]{\PFGstyle $\ph(t,\ )$}%
\psfrag{PhitD}[cr][cr]{\PFGstyle $\ph(t+\Delta,\ )$}%
\psfrag{PhitmD}[cr][cr]{\PFGstyle $\ph(t-\Delta,\ )$}%
\psfrag{TAm1la}[cc][cc]{\PFGstyle $\tfrac{1}{\lambda}(T_A-1)$}%
\psfrag{TEla}[cc][cc]{\PFGstyle $\tfrac{1}{\lambda}T_E$}%
\psfrag{xxA}[Bc][Bc]{\PFGstyle $x$}%
\psfrag{xx}[cl][cl]{\PFGstyle $x$}%
\includegraphics[width=0.5\textwidth]{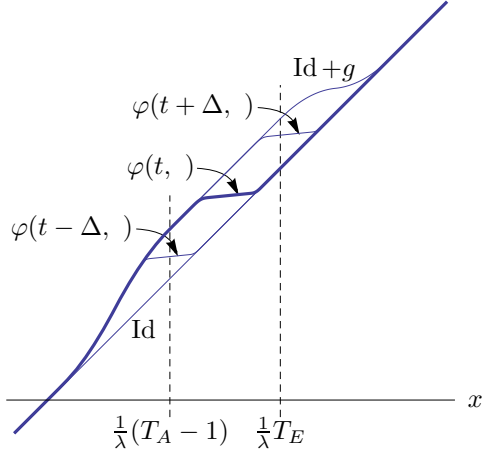}
\end{psfrags}
\caption{The path $\ph(t,\ )$ defined in \nmb!{3.3} connecting $\Id$ to $\Id+g$, 
plotted at $t-\Delta<t<t+\Delta$. 
Between the dashed lines, $g \equiv 1$ is constant. }
\label{fig}
\end{figure}

We consider the function
\begin{equation}\label{eq:f}
\begin{aligned}
 f(z, a, \ep)&=  \max(0, \min(z, a)) \star G_\ep(z) G_\ep(a)\\&=\iint \max(0, \min(z-\ol z, a - \ol a)) G_\ep(\ol z)
G_\ep(\ol a) \ud \ol z \ud \ol a \\
&= \iint \max(0, \min(z - \ep \ol z, a - \ep \ol a)) G_1(\ol z)
G_1(\ol a) \ud \ol z \ud \ol a\\&= \ep f(\frac{z}{\ep}, \frac{a}{\ep}, 1)
\end{aligned}
\end{equation}
where 
$G_\ep(z) = \frac{1}{\ep} G_1(\frac{z}{\ep})$
is a function with $\operatorname{supp}(G_\ep) \subseteq [-\ep,\ep]$ and $\int G_\ep \ud x = 1$. 
Furthermore, let $g :\R\to [0,1]$ be a function with compact support contained in 
$\mathbb R_{>0}$ and $g'>-1$, so that $x+g(x)$ is a diffeomorphism.
For $0<\lambda <1$ and $t \in [0, T]$ let
\[ \ph(t,x) = x + f(t - \lambda x, g(x), \ep(t)) \]
be the path going away from the identity (since $\supp(g)\subset \mathbb R_{>0}$, see also figure 
\ref{fig}). For given $\ep_0>0$, let
\[ \psi(t,x) = x + f(T-t-\lambda x, g(x), \ep_0) \]
the path leading back again. 
The only difference to \cite{Michor2005} is that the parameter $\ep$ may
vary along the path.

We shall need some derivatives of $\ph$ and $f$:
\begin{align*}
\ph_t(t,x) &= f_z(t-\lambda x, g(x), \ep(t)) + \dot \ep(t)
f_\ep(t-\lambda x, g(x), \ep(t))\\
\ph_x(t,x) &=1-\lambda f_z(t-\lambda x, g(x), \ep(t)) 
+f_a(t-\lambda x, g(x), \ep(t))g'(x)\\
f_z(z, a,\ep) &= \int_{-\infty}^z \int_{-\infty}^{a-z} G_\ep(w) G_\ep(w+b)
\ud b \ud w \\
f_a(z, a,\ep) &= \int_{-\infty}^a \int_{-\infty}^{z-a} G_\ep(w) G_\ep(w+b)
\ud b \ud w \\
f_{\ep}(z, a,\ep) &= \frac{1}{\ep} \Big(f(z, a,\ep)-zf_{z}(z, a,\ep)-af_{a}(z, a,\ep)\Big)\\
f_{zz}(z, a,\ep) &= G_\ep(z) \int_{-\infty}^a G_\ep(b) \ud b -
\int_{-\infty}^{z} G_\ep(w) G_\ep(w-(z-a)) \ud w \\
\end{align*}
{\bf Claim 1.} The path $\ph$ followed by $\ps$ still has arbitrarily small length for the $L^2$-metric. 
\newline
We are working with a fixed time interval $[0,2T]$. Thus arbitrarily small length
is equivalent to arbitrarily small energy. The
energy is given by
\begin{equation}\label{eq:e}
\iint \ph_t^2 \ph_x \ud x \ud t=\iint(f_z+\dot \ep f_{\ep})^2(1-\lambda f_z+f_a g') \ud x \ud t
\end{equation}
Looking at the formula for $f_{\ep}$ we see that $\ep f_{\ep}$ is bounded on a domain with bounded $a$. 
Thus $\|\dot \ep f_{\ep}\|_{\infty}\rightarrow 0$ 
can be achieved by choosing $\ep$, such that $|\dot \ep|\leq C\ep^{3/2}$.
We will see  later that this is possible. 
Inspecting $\ph_t(t, x)$
and looking at the formulas for $f_z$ and $f$ we see that for $t - \lambda x <
-\ep(t)$ and for $t - \lambda x - g(x) >2 \ep(t)$ we have $\ph_t(t,x) = 0$. Thus the domain of integration is 
contained in the compact set 
$$[0,T]\times[-\frac{T  + \|g\|_{\infty}+2\|\ep\|_{\infty}}{\lambda},\frac{T  +\|\ep\|_{\infty}}{\lambda}].$$
Therefore, it is enough to show that the $L^{\infty}$-norm of the integrand in \eqref{eq:e} goes to zero 
as $\|\ep\|_{\infty}$ goes to zero. For all terms involving $\dot \ep f_\ep$ this is
true by the above assumption since $(1-\lambda f_z+f_a g')$ and $\ep f_{\ep}$ are bounded.
For the remaining parts $f_z^2(1-\lambda f_z)$ and $f_z^2 f_a g'$ we follow the argumentation of  \cite{Michor2005}.
For $t$ fixed and $\la$ close to $1$, the function $1-\lambda f_z$, when restricted to the support of $f_z$, is bigger than 
$\ep(t)$ only on an
interval of length $O(\ep(t))$. Hence we have
\begin{align*}
\int_0^T\int_{\R} f_z^2(1-\lambda f_z)\ud x \ud t \leq \|f_z\|^2_{\infty}\int_0^T\int_{\R}(1-\lambda f_z)\ud x \ud t
= O(\|\ep\|_{\infty}).
\end{align*}
For the last part, we note that
the support of $f_z^2f_a$ is contained in the set $|g(x) - (t - \la x)| \leq 2\ep$.
Now we define $x_0<x_1$  by $g(x_0) + \la x_0 = T - 2\|\ep\|_{\infty}$ and 
$g(x_1) + \la x_1 = T + 2\|\ep\|_{\infty}$. Then 
\begin{align*}
\int_0^T\int_{\R} f_z^2 f_a g'\ud x \ud t &\leq T\|f_z\|_{\infty}^2\|f_a\|_{\infty} 
\int_{\on{supp}(f_z^2\,f_a)}   g'\ud x
\\&
= T(g(x_1)-g(x_0))\leq 4T \|\ep\|_{\infty}.
\end{align*}
The estimate for $\ps$ is similar and easier.
This proves claim 1.

\noindent
{\bf Claim 2.}
For every $a\in \mathbb R$ and $\de>0$ we may choose $\ep(t)$ with $\|\ep\|_\infty<\de$ such that 
\begin{align*}
\int_0^T\int_{\R} \frac{\ph_{tx}\ph_{xx}}{\ph_x^2} \ud x \ud t +  \int_0^T\int_{\R}
\frac{\psi_{tx}\psi_{xx}}{\psi_x^2} \ud x \ud t = a.
\end{align*}
We will subject $\ep$ and $g$ to several assumptions. First, we partition the
interval $[0,T]$ equidistantly into $0 < T_A  < T_E < T$ 
and the $(t,x)$-domain into two parts, namely $A_1 = ([0, T_A] \cup [T_E, T]) \times \mathbb R$  
and $A_2 = [T_A, T_E]  \times \mathbb R$.
We want $g(x)\equiv 1$ 
on a neighborhood of the interval $[\frac{1}{\lambda}(T_A-1),\frac{1}{\lambda}T_E]$. 
We choose $\ep(t)$ to be constant $\ep(t)
\equiv \ep_0$ on $[0,T_A] \cup [T_E,T]$ and to be symmetric in the sense,
that $\ep(t) = \ep(T - t)$. In addition, we want $\ep(t)$ small enough, such that
$g(x)\equiv 1$ on $[\frac{1}{\lambda}(T_A-1-2\ep(t)),\frac{1}{\lambda}(T_E+\ep(t))].$

On $A_1$ we have
$\ep(t) \equiv \ep_0$. This implies $\psi_{tx}(t,x) =
-\ph_{tx}(T-t,x)$, $\psi_{x}(t,x) = \ph_{x}(T-t,x)$ and $\psi_{xx}(t,x) = \ph_{xx}(T-t,x)$. Hence 
\[ \iint_{A_1} \frac{\ph_{tx}\ph_{xx}}{\ph_x^2} \ud x \ud t + \iint_{A_1}
\frac{\psi_{tx}\psi_{xx}}{\psi_x^2} \ud x \ud t = 0. \]

Let $A_2 = [T_A, T_E]  \times \mathbb R$ be the region, where
$\ep(t)$ is not constant. 
In the interior, where
\[
\begin{array}{rcl}
-\ep(t) &< t - \lambda x <& g(x) + 2\ep(t) \\
t - g(x) - 2\ep(t) &< \lambda x <& t + \ep(t)
\end{array}
\]
we have by assumption $g(x) \equiv 1$. Therefore, one has in this region:
\begin{align*}
\ph_x(t,x) & = -\lambda f_z(t - \lambda x, 1, \ep(t)) \\
\ph_{xx}(t, x) &= \lambda^2 f_{zz}(t - \lambda x, 1, \ep(t)) \\
\ph_{tx}(t, x) &= -\lambda f_{zz}(t - \lambda x, 1, \ep(t)) 
- \lambda f_{\ep z}(t - \lambda x, 1, \ep(t)) \dot \ep(t)
\end{align*}
We divide the integral over $A_2$ into two symmetric parts
\[ \int_{T_A}^{T/2}
\int_{\frac{1}{\lambda}(t-1-2\ep(t))}^{\frac{1}{\lambda}(t+\ep(t))}
\frac{\ph_{tx}\ph_{xx}}{\ph_x^2} \ud x \ud t 
+\int_{T/2}^{T_E}
\int_{\frac{1}{\lambda}(t-1-2\ep(t))}^{\frac{1}{\lambda}(t+\ep(t))}
\frac{\ph_{tx}\ph_{xx}}{\ph_x^2} \ud x \ud t\]
and apply the following variable substitution to the second integral
\begin{align*}
\wt t = T- t,\quad
\wt x = x + \frac{1}{\lambda}(\wt t - t).
\end{align*}
Thus $\wt t - \lambda \wt x = t - \lambda x$. 
Together with $\ep(t)  = \ep(\wt t)$ this implies
$$\ph_x(t, x) = \ph_x(\wt t, \wt x), \quad \ph_{xx}(t, x) = \ph_{xx}(\wt t, \wt x).$$
Since $\dot \ep(t) = - \dot \ep(\wt t)$ changes sign, the term containing $\dot \ep(t)$ 
cancels out and leaves only
$$\ph_{tx}(t, x) + \ph_{tx}(\wt t, \wt x)  = -2 \lambda f_{zz}(t - \lambda x, 1, \ep(t)).$$
A simple calculation shows
that the integration limits transform
\[ \int_{T/2}^{T_E}
\int_{\frac{1}{\lambda}(t-1-2\ep)}^{\frac{1}{\lambda}(t+\ep)}
\frac{\ph_{tx}\ph_{xx}}{\ph_x^2} \ud x \ud t 
= 
\int_{T_A}^{T/2}
\int_{\frac{1}{\lambda}(\wt t-1-2\ep)}^{\frac{1}{\lambda}(\wt t+\ep)}
\frac{\ph_{tx}\ph_{xx}}{\ph_x^2} \ud \wt x \ud \wt t \]
to those of the first integral. 
Therefore, the sum of the integrals gives
\begin{equation*}
\label{eq:intA2} 
\iint_{A_2} \frac{\ph_{tx}\ph_{xx}}{\ph_x^2} \ud x \ud t = -2 \lambda^3
\int_{T_A}^{T/2}
\int_{\frac{1}{\lambda}(t-1-2\ep(t))}^{\frac{1}{\lambda}(t+\ep(t))}
\frac{f_{zz}(t - \lambda x, 1, \ep(t))^2}{\big(1-\lambda f_z(t - \lambda x, 1, \ep(t))\big)^2} \ud x \ud t.
\end{equation*}
From formula \eqref{eq:f} we see:
\begin{align*}
f_{z}(z, a, \ep)& = f_{z}(\frac{z}{\ep}, \frac{a}{\ep}, 1),\qquad
f_{zz}(z, a, \ep) = \frac{1}{\ep} f_{zz}(\frac{z}{\ep},
\frac{a}{\ep}, 1) .
\end{align*}
We can use this to rewrite the above integral:
\begin{align*}
\iint_{A_2} \frac{\ph_{tx}\ph_{xx}}{\ph_x^2} \ud x \ud t
&=-2 \lambda^3\int_{T_A}^{T/2}
\int_{\frac{1}{\lambda}(t-1-2\ep(t))}^{\frac{1}{\lambda}(t+\ep(t))}
\frac{f_{zz}(t - \lambda x, 1, \ep(t))^2}{\big(1-\lambda f_z(t - \lambda x, 1, \ep(t))\big)^2} \ud x \ud t\\
&=-2 \lambda^2
\int_{T_A}^{T/2}
\int_{-\ep(t)}^{2\ep(t)+1}
\frac{f_{zz}(z, 1, \ep(t))^2}{\big(1-\lambda f_z(z, 1, \ep(t))\big)^2} \ud z \ud t\\
&=-2 \lambda^2
\int_{T_A}^{T/2}
\int_{-\ep(t)}^{2\ep(t)+1}\frac{1}{\ep(t)^2}
\frac{f_{zz}(\tfrac{z}{\ep(t)},\tfrac{1}{\ep(t)},1)^2}
{\big(1-\lambda f_z(\tfrac{z}{\ep(t)}, \tfrac{1}{\ep(t)}, 1)\big)^2} \ud z \ud t\\
&=-2 \lambda^2
\int_{T_A}^{T/2}
\int_{-1}^{2+\frac{1}{\ep(t)}}\frac{1}{\ep(t)}
\frac{f_{zz}(z,\tfrac{1}{\ep(t)},1)^2}
{\big(1-\lambda f_z(z, \tfrac{1}{\ep(t)}, 1)\big)^2} \ud z \ud t
\end{align*}
Looking at the formula for $f_{zz}$ 
\[ f_{zz}(z,\tfrac{1}{\ep},1) = G_1(z) - \int_{-\infty}^{z} G_1(w) G_1(w - (z - \tfrac{1}{\ep})) \ud w \]
we see that $f_{zz}(z,\tfrac{1}{\ep},1)$ is non-zero only on the intervals 
$|z| < 1$ and $|z - \tfrac{1}{\ep}| < 2$. 
For small $\ep$, these are two disjoint regions. Therefore, the above integral equals
\begin{align*}
\iint_{A_2} \frac{\ph_{tx}\ph_{xx}}{\ph_x^2} \ud x \ud t =&-2 \lambda^2
\int_{T_A}^{T/2} \frac{1}{\ep(t)}
\int_{-1}^1 \frac{f_{zz}(z,\tfrac{1}{\ep(t)},1)^2}
{\big(1-\lambda f_z(z, \tfrac{1}{\ep(t)}, 1)\big)^2} \ud z \ud t- \\ 
&-2 \lambda^2
\int_{T_A}^{T/2} \frac{1}{\ep(t)}
\int_{-2}^2 \frac{f_{zz}(z+\tfrac{1}{\ep(t)},\tfrac{1}{\ep(t)},1)^2}
{\big(1-\lambda f_z(z+\tfrac{1}{\ep(t)}, \tfrac{1}{\ep(t)}, 1)\big)^2} \ud z \ud t 
\end{align*}
For $z$ bounded and sufficiently small $\ep(t)$, 
the functions under the integral do not depend on $\ep(t)$ any more
as can be seen from the definitions of $f_z$ and $f_{zz}$.  
Thus
$$I = \lambda^2 \int_{-1}^1 \frac{f_{zz}(z,\tfrac{1}{\ep(t)},1)^2}
{\big(1-\lambda f_z(z, \tfrac{1}{\ep(t)}, 1)\big)^2} \ud z 
+\lambda^2 \int_{-2}^2 \frac{f_{zz}(z+\tfrac{1}{\ep(t)},\tfrac{1}{\ep(t)},1)^2}
{\big(1-\lambda f_z(z+\tfrac{1}{\ep(t)}, \tfrac{1}{\ep(t)}, 1)\big)^2} \ud z,
$$
is independent of $t$ and we have
\begin{align*}
\iint_{A_2} \frac{\ph_{tx}\ph_{xx}}{\ph_x^2} \ud x \ud t  =
- I \int_{T_A}^{T_E}\frac{1}{\ep(t)}\ud t.
\end{align*}

The same calculations can be repeated for the return path $\psi$, where
$\ep \equiv \ep_0$ is constant in time:
\[ \iint_{A_2} \frac{\psi_{tx}\psi_{xx}}{\psi_x^2} \ud x \ud t =
 I \int_{T_A}^{T_E}\frac{1}{\ep_0}\ud t. \]
Note that the sign is positive now, which comes from the $t$-derivative. Putting
everything together gives us
\begin{equation*}
a = \iint \frac{\ph_{tx}\ph_{xx}}{\ph_x^2} + \frac{\psi_{tx}\psi_{xx}}{\psi_x^2}
\ud x \ud t = I\int_{T_A}^{T_E}\Big( \frac{1}{\ep_0}-\frac{1}{\ep(t)}\Big)\ud t 
\end{equation*}
Let $\ep(t)=\ep_0+\ep_1\ep_0^{3/2}b(t)$ where
$b(t)$ is a bump function with height $1$ and $\ep_1$  is a small constant. 
Note that  $\ep(t)$ satisfies  $|\dot \ep|\leq \|\dot b\|_\infty\,\ep_1 \ep_0^{3/2}$.
Choosing $\ep_0$ and $\ep_1$ small independently we may produce any $a\in \mathbb R$.
\qed\end{demo}

\begin{demo}{Proof of Theorem \nmb!{3.1}}
Let $(\ph,a)\in \mathbb R\x_c\on{Diff}_{\mathcal S}(\mathbb R)$. By proposition \nmb!{3.2} we get a 
smooth family $\ph(\de, t, x)$ for $\de>0$ and $t\in [0,1]$ such that 
$\ph(\de,t,\quad)\in\on{Diff}_{\mathcal S}(\mathbb R)$,
$\ph(\de,0,\quad)=\on{Id}_{\mathbb R}$,
$\ph(\de,1,\quad)=\ph$, and such that
the length of $t\mapsto \ph(\de,t,\quad)$ is $<\de$.

Using \nmb!{2.4} consider the horizontal lift 
$(\ph(\de,t,\quad),a(\de,t))\in\on{Diff}_{\mathcal S}(\mathbb R)$ of this family
which connects $\binom{\Id}{0}$ with $\binom{\ph}{a(\de,1)}$ for each $\de>0$
and has length $<\de$.
But one can see from the proof of lemma \nmb!{3.3} that $a(\de,1)$ becomes unbounded for 
$\de\to 0$. 

Using lemma \nmb!{3.3} we can find a horizontal path 
$t\mapsto \binom{\ps(\de,t,\quad)}{b(\de,t)}$ for $t\in[0,1]$ 
in the Virasoro group of length $<\de$ connecting 
$\binom{\Id}{0}$ with $\binom{\Id}{a-a(\de,1)}$. Then the curve 
$t\mapsto \binom{\ps(\de,t,\quad)}{b(\de,t)}.\binom{\ph}{a(\de,1)}
=\binom{\ps(\de,t)\o\ph}{b(\de,t)+a(\de,1)+c(\ps(\de,t),\ph)}$
connects
$\binom{\ph}{a(\de,1)}=\binom{\Id}{0}.\binom{\ph}{a(\de,1)}$ with 
$\binom{\ph}{a}=\binom{\Id}{a-a(\de,1)}.\binom{\ph}{a(\de,1)}$ and it has length $<\de$.
\qed\end{demo}

\bibliographystyle{plain}

\begin{thebibliography}{10}

\bibitem{Arnold66} 
V.I. Arnold:
\newblock Sur la g\'eometrie diff\'erentielle des groupes de Lie de  
dimension infinie et ses applications \`a l'hydrodynamique des  
fluides parfaits. 
\newblock {\em Ann. Inst. Fourier} 16 (1966), 319--361. 

\bibitem{Constantin2003}
Adrian Constantin and Boris Kolev.
\newblock Geodesic flow on the diffeomorphism group of the circle.
\newblock \emph{Comment. Math. Helv.}, 78\penalty0 (4):\penalty0 787--804,
  2003.

\bibitem{CKKT07}
A.~Constantin, T.~Kappeler, B.~Kolev, and P.~Topalov.
\newblock On geodesic exponential maps of the Virasoro group.
\newblock  \emph{Ann. Global Anal. Geom.} 31 (2007), 155�180.


\bibitem{Fuks86} 
D.~Fuks,  
\newblock {\em Cohomology {o}f {i}nfinite {d}imensional {L}ie algebras}.
\newblock Nauka, Moscow, 1984 (Russian).
\newblock Transl. English, Contemporary Soviet Mathematics.
Consultants Bureau (Plenum Press), New York, 1986.

\bibitem{KM97}
A. Kriegl, P.W. Michor:
\newblock {\em The Convenient Setting for Global Analysis.}   
\newblock Surveys and Monographs 53, AMS, Providence 1997.

\bibitem{MGeomEvol06}
P.~W. Michor. 
\newblock Some Geometric Evolution Equations Arising as Geodesic Equations 
on Groups of Diffeomorphism, Including the Hamiltonian Approach. 
\newblock In: Phase Space Analysis of Partial Differential Equations. 
Bove, Antonio; Colombini, Ferruccio; Santo, Daniele Del (Eds.).
\newblock Progress in Non Linear Differential Equations and 
Their Applications, Vol. 69.  
\newblock Birkh\"auser Boston, 2006.

\bibitem{MM06}
Peter~W. Michor and David Mumford.
\newblock Riemannian geometries on spaces of plane curves.
\newblock \emph{J. Eur. Math. Soc. (JEMS) 8 (2006), 1-48}, 2006.

\bibitem{Michor2005}
Peter~W. Michor and David Mumford.
\newblock Vanishing geodesic distance on spaces of submanifolds and
  diffeomorphisms.
\newblock {\em Documenta Math.}, 10:217--245, 2005.

\bibitem{MR98} 
P.W. Michor, T. Ratiu:
\newblock On the geometry of the Virasoro-Bott group.  
\newblock {\em  J. Lie Theory} 8, 2 (1998), 293-309.    

\bibitem{Misiolek97} 
G.~Misio\l ek. 
\newblock Conjugate points in the Bott-Virasoro group and the KdV equation.
\newblock{\em Proc. Amer. Math. Soc.}  
125 (1997), 935--940.

\bibitem{OvsienkoKhesin87}
V.Y. Ovsienko and B.A. Khesin.
\newblock Korteweg--de~Vries superequations as an Euler equation.
\newblock {\em Funct. Anal. Appl.}, 21 (1987), 329--331.

\bibitem{Schilling2006}
Ren{\'e} L. Schilling:
\newblock {\em Measures, Integrals and Martingales.}
\newblock Cambridge University Press, New York 2005.

\end{thebibliography}

\end{document}